\documentclass[amsfonts,11pt]{amsart}
\RequirePackage{amsmath,amssymb,amsthm, amscd, comment,mathtools}
\usepackage[margin=1.2in]{geometry}
\usepackage{amsmath}
\usepackage{amsthm}
\usepackage{amssymb}
\usepackage{amscd}
\usepackage{mathrsfs}
\usepackage{graphicx}
\usepackage{rotating, graphicx}
\usepackage[all,cmtip]{xy}
\usepackage{enumerate}
\usepackage{tabularx}
\usepackage{tikz-cd}
\usepackage[thinlines]{easytable}
\usepackage{multirow}
\usepackage{mathtools}
\usepackage{mathrsfs}
\usepackage{changepage} 
\usepackage{paralist}

\usepackage{kantlipsum}
\usepackage{mathabx}
\usepackage{enumitem}
\usepackage{times}
\usepackage{subfig}
\usepackage{color}
\usepackage[pagebackref,breaklinks,colorlinks,linkcolor=red,anchorcolor=red,citecolor=blue]{hyperref}

\calclayout

\newcommand{\Gm}{\mathbb{G}_m}

\newtheorem{proposition}[subsection]{Proposition}
\newtheorem{corollary}[subsection]{Corollary}
\newtheorem{theorem}[subsection]{Theorem}
\newtheorem{lemma}[subsection]{Lemma}
\newtheorem{conjecture}[subsection]{Conjecture}
\theoremstyle{definition}

\theoremstyle{remark}

\newtheorem{remark}[subsection]{Remark}

\numberwithin{equation}{subsection}

\DeclareMathAlphabet{\mathbbold}{U}{bbold}{m}{n}

\title[]{Motivic absolute purity away from the characteristic}

\author{Fangzhou Jin}
\address{School of Mathematical Sciences\\
Key Laboratory of Intelligent Computing and Applications (Ministry of Education)\\
Tongji University\\
Siping Road 1239\\
200092 Shanghai\\
China}
\email{\href{mailto:fangzhoujin@tongji.edu.cn}{fangzhoujin@tongji.edu.cn}}
\urladdr{\url{https://fangzhoujin.github.io/}}

\date{\number\day-\number\month-\number\year}

\begin{document}

\maketitle

\begin{abstract}

We show that the motivic sphere spectrum satisfies the absolute purity property after inverting the residual characteristics, proving a conjecture of D\'eglise in that case. Consequently the same property holds for every cellular spectrum, such as $\mathbf{MGL}$ or $\mathbf{H}\mathbb{Z}$.

\end{abstract}

\setcounter{tocdepth}{1}
\tableofcontents

\noindent
\section{Introduction}

\subsection{Absolute purity for \'etale sheaves}
Grothendieck's \emph{absolute (cohomological) purity conjecture} (\cite[Expos\'e I 3.1.4]{SGA5}) is the following statement: if $i:Z\to X$ is a closed immersion between noetherian regular schemes of pure codimension $c$, $n\in\mathcal{O}(X)^*$ and $\Lambda=\mathbb{Z}/n\mathbb{Z}$, then the \'etale cohomology sheaf supported in $Z$ with values in $\Lambda$ can be computed as
\begin{gather}
\begin{align}
\mathcal{H}^q_Z(X_{\textrm{\'et}},\Lambda)=
  \begin{cases} 
   i_*\Lambda_Z(-c) & \text{if } q=2c, \\
   0 & \text{else}.
  \end{cases}
\end{align}
\end{gather}
In other words, we have an isomorphism $i^!\Lambda_X\simeq\Lambda_Z(-c)[-2c]$. This conjecture is eventually settled in full generality by Gabber (see \cite{Fuj}), using Atiyah-Hirzebruch spectral sequence for \'etale K-theory developed by Thomason (\cite{Tho}). Later Gabber found a second proof based on log alterations and a computation for semistable models (\cite[Exp. XVI]{ILO}).


\subsection{Purity in motivic homotopy}
\label{recall_pur_trans}
The absolute purity property has the following natural generalization to motivic homotopy.
By \cite[4.3.7]{DJK}, for any smoothable lci morphism $f\colon X\to S$ with cotangent complex $L_f$ and any motivic spectrum $E\in\mathbf{SH}(S)$, there is a \emph{purity transformation}
\begin{equation}
\label{eq:pur_trans}
f^*E\otimes \operatorname{Th}_X(L_f)\to f^!E.
\end{equation}
We say that the object $E$ is \textbf{pure relatively to $f$} or \textbf{$f$-pure} if the map~\eqref{eq:pur_trans} is an isomorphism. 

An \emph{absolute motivic spectrum} $E$ (that is, the data of a spectrum $E_X\in\mathbf{SH}(X)$ for every scheme $X$, such that $f^*E_X\simeq E_Y$ for every morphism $f:Y\to X$) satisfies \textbf{absolute purity} if for every closed immersion $i:Z\to X$ between regular schemes, $E_X$ is $i$-pure (see \cite[Appendix C]{DFJK}). This property has important applications, including proving the constructibility of the six functors and Grothendieck-Verdier local duality (see \cite{DFJK}), as well as in the study of the coniveau spectral sequence (see \cite{Jin16}).
In \cite[Conjecture B]{Deg11}, D\'eglise made the following conjecture:
\begin{conjecture}
\label{conj:sphereap}
The motivic sphere spectrum $\mathbf{1}$ satisfies absolute purity. 
\end{conjecture}
By \cite[Remark 1.3.5]{Deg11}, Conjecture~\ref{conj:sphereap} implies that every cellular motivic spectrum also satisfies absolute purity. 
The equal characteristic case follows from Popescu's theorem (\cite{Spi}), see \cite[Appendix C]{DFJK}, so the mixed characteristic case is where lies the real difficulty.
Below is a list of known cases: 
\begin{itemize}
\item 
Cisinski-D\'eglise: (homotopy) algebraic K-theory spectrum $\mathbf{KGL}$ and rational motivic Eilenberg-Mac Lane spectrum $\mathbf{H}\mathbb{Q}$ (\cite[Thm. 13.6.3, 14.4.1]{CD})
\item
D\'eglise: rational algebraic cobordism spectrum $\mathbf{MGL}_{\mathbb{Q}}$ (\cite[Example 1.3.4 (4)]{Deg11})
\item
D\'eglise-Fasel-Jin-Khan: rational motivic sphere spectrum $\mathbf{1}_\mathbb{Q}$ (\cite[Cor. C]{DFJK})
\item
Calm\`es-Harpaz-Nardin: Hermitian K-theory spectrum $\mathbf{KQ}$ (\cite[Thm. 8.4.2]{CHN})
\item
R\"ondigs-Spitzweck: $\mathbf{1}[p^{-1}]$ is pure with respect to the inclusion of a closed point in the spectrum of a discrete valuation ring of mixed characteristic $(0,p)$ (\cite[Thm. 7.8]{RS})
\item
Frankland-Spitzweck: motivic Eilenberg-Mac Lane spectrum $\mathbf{H}\mathbb{Z}$ (as defined in \cite{Spit}) is pure with respect to the inclusion of a closed point in a Dedekind scheme (\cite[Cor. 8.4]{FS})
\item
R\"ondigs: $\mathbf{1}[\eta^{-1}]$ is pure with respect to the inclusion of a closed point in a Dedekind scheme (\cite[Cor. 5.10]{Ron})
\end{itemize}

The main result of this paper is to prove Conjecture~\ref{conj:sphereap} after inverting the characteristics:
\begin{theorem}
\label{thm:pur1}
Let $i:Z\to X$ be a closed immersion between regular noetherian schemes, and denote by $\operatorname{char}(X)$ the set of prime numbers not invertible on $X$. Then the motivic sphere spectrum $\mathbf{1}_X[\operatorname{char}(X)]^{-1}$ is $i$-pure.
\end{theorem}
Consequently, by \cite[Remark 1.3.5]{Deg11}, we deduce the following:
\begin{corollary}
Let $\mathcal{P}$ be a collection of prime numbers. For any cellular motivic spectrum $E$ such as $\mathbf{MGL}$ or $\mathbf{H}\mathbb{Z}$, the $\mathcal{P}$-inverted spectrum $E[\mathcal{P}^{-1}]$ satisfies absolute purity when restricted to schemes $X$ such that $\operatorname{char}(X)\subset\mathcal{P}$.
\end{corollary}


%

\subsection{Proof strategy}
In \cite{DFJK}, absolute purity for the rational motivic sphere spectrum $\mathbf{1}_\mathbb{Q}$ is obtained by decomposing it into the direct sum of the plus and minus parts $\mathbf{1}_{+,\mathbb{Q}}$ and $\mathbf{1}_{-,\mathbb{Q}}$. The decomposition exists as long as $2$ is inverted, and the $2$-inverted minus part $\mathbf{1}_{-}[2^{-1}]$ is concentrated on the characteristic $0$ fiber by Bachmann's theorem (\cite[Theorem 35]{Bac}), and therefore its absolute purity is reduced to the equal characteristic case. For the plus part, Cisinski-D\'eglise identify its rationalization $\mathbf{1}_{+,\mathbb{Q}}$ with $\mathbf{H}\mathbb{Q}$ and deduce its purity from d\'evissage for algebraic K-theory (\cite{CD}). Without rationalization, unfortunately, not much is known about the plus part.

Our Theorem~\ref{thm:pur1} does not even require inverting $2$ if $2\notin\operatorname{char}(X)$, and the strategy is a refinement of the approach above. We start by reducing the problem to purity of $\mathbf{1}[p^{-1}]$ between regular schemes of finite type over $\mathbb{Z}_{(p)}$ in Proposition~\ref{prop:Zpfin}, by applying a theorem of Popescu (\cite[Thm. 3.8]{Pop}, which is refined in Lemma~\ref{lm:DFJKC8} below). We further reduce to the $\eta$-quotient sphere $(\mathbf{1}/\eta)[p^{-1}]$ and the $\eta$-inverted sphere $\mathbf{1}[\eta^{-1}]$ respectively. Since the purity for the $(2,\eta)$-inverted sphere $\mathbf{1}[2^{-1},\eta^{-1}]$ is known by \cite[Thm. 3.6]{DFJK}, we are reduced to prove purity for the $2$-local $\eta$-inverted sphere $\mathbf{1}[\eta^{-1}]_{(2)}$ when $p$ is odd.

To prove purity for the $\eta$-quotient sphere $(\mathbf{1}/\eta)[p^{-1}]$, we first reduce it to purity for $\mathbf{MGL}[p^{-1}]$ in Proposition~\ref{prop:HF2pur}: indeed, work of Mantovani (\cite{Mon}) identifies $\mathbf{1}/\eta$ with its $\mathbf{MGL}$-cobar resolution. The case of $\mathbf{MGL}[p^{-1}]$ is proven in Corollary~\ref{cor:MGLpur}, where the proof follows the same lines as Gabber's second proof in \cite[Exp. XVI]{ILO}: using log alterations (\cite[Thm. 4.3.1]{Tem}), one reduces to a log smooth situation (Lemma~\ref{lm:Gabber362}), and eventually to an explicit computation in the semistable case (Lemma~\ref{lm:chartpur}). Our computation in this case differs from Gabber's approach, as we do not want to compute \'etale locally; instead, we use torus actions to circumvent the difficulty, see Lemma~\ref{lm:eqret} and Pham's work \cite[Lemma 2.23]{Pha}. These computations are similar to Ayoub's work in \cite[Thm. 3.3.10 and 3.3.24]{Ayo}.
In Lemma~\ref{lm:MGLpurdes} we show that after inverting the degree, purity for $\mathbf{MGL}$ desends through alterations, which allows to conclude. The argument relies on orientability of $\mathbf{MGL}$ as well as a degree computation, which also works for $\mathbf{KGL}$ or $\mathbf{H}\mathbb{Z}$ but unfortunately not for the sphere spectrum. Note that recent work of Bachmann-Hopkins (\cite{BH}) implies that purity for $\mathbf{MGL}$ with integral coefficients is equivalent to purity for $\mathbf{H}\mathbb{Z}$, and consequently they hold over a Dedekind scheme, see Lemma~\ref{lm:purMGLHZ} and Corollary~\ref{cor:MGLDpur}.

It remains to prove purity for the $2$-local $\eta$-inverted sphere $\mathbf{1}[\eta^{-1}]_{(2)}$ over $\mathbb{Z}[1/2]$-schemes. By Bachmann's fiber sequence (\cite[Theorem 4.12]{Bac22}, see~\eqref{eq:Bac412} below), we only need to prove purity for $\mathbf{kw}$, the connective cover of the periodic Witt spectrum in~\ref{num:Bacsp} (see also \cite{Ron}); the last three isomorphisms in~\eqref{eq:KWrel} reduce it to $\mathbf{K}^W$. On the other hand, the first isomorphism in~\eqref{eq:KWrel} together with purity for $\mathbf{H}\mathbb{F}_2$ (proven in Proposition~\ref{prop:HF2pur} and Corollary~\ref{cor:MGLpur}) shows purity for $\mathbf{K}^W/\eta$. To lift purity to $\mathbf{K}^W$, it suffices to show that $\mathbf{K}^W$ is $\eta$-complete. This is not true over general base schemes, but 
we can reduce to the case of a carefully chosen finite \'etale $\mathbb{Z}_{(p)}$-algebra $D=\mathbb{Z}_{(p)}[u]/(u^2+p-1)$ which is a Nisnevich neighborhood of the closed point of $\mathbb{Z}_{(p)}$, with no orderable residue fields. By \cite[Prop. 10.19]{BEM}, we are reduced to a vanishing statement for the cdh cohomology of the powers of the fundamental ideal $\underline{I}_{\operatorname{cdh}}^m$, which is treated in Lemma~\ref{lm:Ijfin}. Consequently we prove $\eta$-completeness of $\mathbf{K}^W$ for $D$-schemes of finite type in Lemma~\ref{lm:KWcomp}. Finally we prove purity for $\mathbf{1}[\eta^{-1}]_{(2)}$ in Proposition~\ref{prop:1eta2pur}.

To summarize, Theorem~\ref{thm:pur1} is deduced from Proposition~\ref{prop:Zpfin}, Lemma~\ref{lm:1pred}, Proposition~\ref{prop:HF2pur}, Corollary~\ref{cor:MGLpur} and Proposition~\ref{prop:1eta2pur} below, each of which can be checked individually.

\subsection{Acknowledgments}

The author would like to thank Fr\'ed\'eric D\'eglise, Jean Fasel, Niels Feld and Adeel Khan for previous collaborations and discussions that inspired this work. He is supported by 
the National Natural Science Foundation of China Grant Nr.12471014.

\subsection{AI disclaimer}

Commercially available versions of large language models (LLMs) have been used for searching the literature and helping the author understand certain technical points, such as Gabber's proof for absolute purity. The text of this article is entirely written by human, and the author assumes full responsibility for the mathematical statements, proofs, and references.

\section{Notations and preliminaries}

\subsection{}
Throughout the paper, all schemes are assumed noetherian. Given a scheme $X$, denote by $\operatorname{char}(X)$ the set of prime numbers not invertible on $X$. A separated morphism of finite type is \emph{smoothable} if it admits a global factorization through a closed immersion and a smooth morphism; such a morphism is \emph{smoothable lci} if the corresponding closed immersion is a regular closed immersion.

\subsection{}
For $f:X\to S$ a morphism of schemes and $E\in\mathbf{SH}(S)$, we denote $E_{|X}=f^*\mathbb{E}\in\mathbf{SH}(X)$.

\subsection{}
The following lemma follows from transitivity of the purity transformation~\eqref{eq:pur_trans} (see \cite[Prop. 4.3.10]{DJK}):
\begin{lemma}
\label{lm:DJK4310}
Let $h:Y\to X$ and $g:X\to S$ be two smoothable lci morphism, such that the composition $gh$ is smoothable lci. If $E\in\mathbf{SH}(S)$ is both $g$-pure and $gh$-pure, then $g^*E$ is $h$-pure.
\end{lemma}

\subsection{}
For $f:X\to S$ a smoothable lci morphism and $E\in\mathbf{SH}(S)$ a motivic spectrum, denote by $f^\triangle E$ the cofiber of the natural transformation~\eqref{eq:pur_trans}:
\begin{align}
f^\triangle E=\operatorname{cofib}(f^*E\otimes \operatorname{Th}_X(L_f)\xrightarrow{\eqref{eq:pur_trans}} f^!E).
\end{align}
One can define $f^\triangle$ as an exact functor, see \cite[Remark 6.1.4]{JY}.
In the case of a closed immersion between regular schemes, the object $f^\triangle E$ measures the defect of absolute purity.

\subsection{}
Denote by $\eta:\mathbf{1}(1)[1]\to \mathbf{1}$ the motivic Hopf map, obtained by delooping the canonical morphism $\mathbb{A}^2\setminus\{0\}\to\mathbb{P}^1$. If $E$ is a motivic spectrum, denote by $E/\eta=\operatorname{cofib}(E(1)[1]\xrightarrow{\eta}E)$ the $\eta$-quotient, and $E[\eta^{-1}]=\operatorname{colim}\left(E\xrightarrow{\eta}E(-1)[-1]\xrightarrow{\eta}E(-2)[-2]\cdots\right)$ the $\eta$-inversion.

\subsection{}
Let $\mathcal{P}$ be a set of prime numbers, and let $M_{\mathcal{P}}\subset \mathbb{N}_{>0}$ be the multiplicative subset generated by $\mathcal{P}$. If $E$ is a motivic spectrum, the localization $E[\mathcal{P}^{-1}]$ is defined as the colimit $E[\mathcal{P}^{-1}]=\operatorname{colim}_{m\in M_{\mathcal{P}}}E$, where the transition map from the copy indexed by $m$ to the copy indexed by $mn$ is the multiplication by $n$. In particular, if $p$ is a prime number, let $\mathcal{P}_p$ be the set of all primes except $p$, and denote by $E_{(p)}=E[\mathcal{P}_p^{-1}]$ the $p$-localization.

\subsection{}
The following lemma shows that purity is a Nisnevich local condition:
\begin{lemma}
\label{lm:purbc}
Let $f:X\to S$ be a smoothable lci morphism and let $E\in\mathbf{SH}(S)$.
\begin{enumerate}
\item
If $E$ is $f$-pure, then for any smooth morphism $Y\to S$, $E_{|Y}$ is pure relatively to the base change $f_Y:X\times_SY\to Y$ of $f$.
\item
Let $(U_i\to S)_{i\in I}$ be a Nisnevich cover of $S$. Then $E$ is $f$-pure if and only if for any $i\in I$, $E_{|U_i}$ is pure relatively to the base change $f_i:X\times_SU_i\to U_i$ of $f$.
\end{enumerate}
\end{lemma}
The spectra $\mathbf{MGL}$ and $\mathbf{H}\mathbb{Z}$ are cellular with integral coefficients over any scheme (\cite[Thm. 6.4]{DI}, \cite[Theorem 7.2]{BH}).

%

\subsection{}
An object of $\mathbf{SH}(S)$ is \textbf{cellular} if it lies in the smallest localizing stable subcategory which contains all the bigraded spheres $\mathbf{1}(q)[p]$. Cellular objects are obtained from the spheres $\mathbf{1}(q)$ by taking shifts, retracts and colimits, and since purity is preserved by colimits by \cite[Lemma 7.9]{RS}, we deduce the following:
\begin{lemma}
\label{lm:celpur}
Let $f:X\to S$ be a smoothable lci morphism and let $E,F\in\mathbf{SH}(S)$ be two objects. If $E$ is $f$-pure and if $F$ is cellular, then $E\otimes F$ is $f$-pure.
\end{lemma}

\subsection{}
The following lemma states conservativity of $\mathbf{SH}$ with respect to purely transcendental field extensions:
\begin{lemma}
\label{lm:extcons}
Let $k$ be a field and denote by $k(t)$ the field of rational fractions with one variable. If $E\in\mathbf{SH}(k)$ is such that $E_{|k(t)}=0$, then $E=0$. 
\end{lemma}
\proof
Let $\mathcal{O}$ be the local ring of $\mathbb{A}^1_k$ at $0$, and denote by $i$ the inclusion of the closed point in $\mathcal{O}$. Since $E_{|k(t)}=0$, the localization sequence induces an isomorphism $i^!E_{|\mathcal{O}}\simeq E$. On the other hand, by relative purity, this map is identified with the map $E(-1)[-2]\to E$ induced by the Euler class of the normal bundle $N_i$ of $i$. Since $N_i$ is a trivial bundle, its Euler class is zero, and therefore $E=0$.
\endproof

\subsection{}
Recall the following lemma from \cite[Cor. 14]{Bac}. See also \cite[Prop. B.20]{DFJK} for a variant.
\begin{lemma}
\label{lm:ptcons}
Let $X$ be a noetherian scheme. As $x$ ranges over the points of $X$, the family of functors
\begin{align}
i_x^*:\mathbf{SH}(X)\to\mathbf{SH}(\operatorname{Spec}\kappa(x))
\end{align}
is jointly conservative.
\end{lemma}
Indeed, the finite dimensionality assumption in loc. cit. can be dropped: by \cite[Prop. 4.3.9]{CD} we can reduce the statement for $X$ to all its local rings, and since every noetherian local ring has finite dimension (\cite[Tag 00KQ]{Stack}), one can use induction on the dimension of local rings and proceed as in the proof of \cite[Cor. 14]{Bac}.


\section{First reductions}

\subsection{}
Recall Popescu's theorem (\cite[Thm. 3.8]{Pop}): any regular local ring of mixed characteristic is a filtered colimit of regular local rings essentially of finite type over $\mathbb{Z}$. Recall from \cite[Def. C.3]{DFJK} that a \emph{regular closed pair} $(X,Z)$ is a regular closed immersion $i:Z\to X$ between schemes; a \emph{transverse morphism} of regular closed pairs $f:(Y,W)\to (X,Z)$ is a Tor-independent cartesian square
\begin{align}
\begin{split}
  \xymatrix@=10pt{
    W \ar[r]^-{} \ar[d]_-{} & Y \ar[d]^-{} \\
    Z \ar[r]_-{} & X.
  }
\end{split}
\end{align} 
The following is a mixed characteristic version of \cite[Prop. C.8]{DFJK}:
\begin{lemma}
\label{lm:DFJKC8}
Let $A$ be a regular local ring of mixed characteristic $(0,p)$, with maximal ideal $\mathfrak{m}$. Let $I\subset A$ be an ideal such that $A/I$ is regular. Then the pair $(\operatorname{Spec}A,\operatorname{Spec}(A/I))$ is the cofiltered limit of regular closed pairs $(\operatorname{Spec}A_\alpha,\operatorname{Spec}(A_\alpha/I_\alpha))$ with transverse transition morphisms, such that each $A_\alpha$ is a regular local ring essentially of finite type over $\mathbb{Z}_{(p)}$.
\end{lemma}
\proof
By Popescu's theorem, write $A=\operatorname{colim}A_\alpha$ with each $A_\alpha$ regular local and essentially of finite type over $\mathbb{Z}$. Replacing each $A_\alpha$ by its localization at the inverse image of $\mathfrak{m}$, we still obtain regular local rings with local transition maps, essentielly of finite type over $\mathbb{Z}_{(p)}$, and the colimit is still $A$ by universal property of localization; in addition, each ring map $A_\alpha\to A$ is local.

Let $c=\dim A-\dim(A/I)$. Choose $t_1,\cdots,t_c\in I$ as part of a regular system of parameters of $A$ such that $I=(t_1,\cdots,t_c)$ (\cite[Tag 00NR]{Stack}). Lift the $t_j$'s to some stage $A_0$, and for any $\alpha\geqslant0$ denote by $t_{j,\alpha}$ the image in $A_\alpha$, and denote $I_\alpha=(t_{1,\alpha},\cdots,t_{c,\alpha})$. By \cite[Lemma C.7]{DFJK}, we know that the $t_{j,\alpha}$'s form a regular system of parameters of $A_\alpha$, and $A_\alpha/I_\alpha$ is regular.

For any $\beta\geqslant\alpha$ we have $I_\beta=I_\alpha A_\beta$, so the morphism of pairs 
\begin{align}
(\operatorname{Spec}A_\beta,\operatorname{Spec}(A_\beta/I_\beta))\to (\operatorname{Spec}A_\alpha,\operatorname{Spec}(A_\alpha/I_\alpha))
\end{align} 
gives rise to a Cartesian diagram. Furthermore, tensoring the Koszul resolution of $A_\alpha/I_\alpha$ with $A_\beta$ gives the Koszul resolution of $A_\beta/I_\beta$, which shows that the square is Tor-independent, hence transverse (\cite[Def. C.3]{DFJK}). Since $A/I=\operatorname{colim}_\alpha(A_\alpha/I_\alpha)$, the result follows.
\endproof




\begin{proposition}
\label{prop:Zpfin}
Assume that for every prime number $p$, the motivic sphere $\mathbf{1}[p^{-1}]$ is pure relatively to every closed immersion between regular schemes of finite type over $\mathbb{Z}_{(p)}$. Then Theorem~\ref{thm:pur1} holds.
\end{proposition}
\proof
Let $i:Z\to X$ be a closed immersion between regular noetherian schemes. Arguing as in the second proof of \cite[Thm. C.1]{DFJK},
it suffices to show that for every point $z\in Z$ with residual characteristic $p$, the sphere $\mathbf{1}[p]^{-1}$ is pure relatively to the closed immersion $i_z:\operatorname{Spec}\mathcal{O}_{Z,z}\to\operatorname{Spec}\mathcal{O}_{X,z}$. 
If $\mathcal{O}_{X,z}$ has equal characteristic, the result follows from \cite[Thm. C.1]{DFJK}. Therefore we may assume $\mathcal{O}_{X,z}$ has mixed characteristic $(0,p)$. By Lemma~\ref{lm:DFJKC8} and our assumption, the pair $(\operatorname{Spec}\mathcal{O}_{X,z},\operatorname{Spec}\mathcal{O}_{Z,z})$ is the cofiltered limit of pairs $(X_\alpha,Z_\alpha)$ with affine and transverse transition morphisms, and the sphere $\mathbf{1}[p^{-1}]$ is pure relatively to each pair $(X_\alpha,Z_\alpha)$. Then the result follows from \cite[Lemma C.4]{DFJK}.
\endproof

%


%

\begin{lemma}
\label{lm:1pred}
Let $p$ be a prime number and let 
$f$ be a smoothable lci morphism.
Assume that the spectrum $(\mathbf{1}/\eta)[p^{-1}]$ is $f$-pure. If $p$ is odd, assume in addition that the spectrum $\mathbf{1}[\eta^{-1}]_{(2)}$ is $f$-pure. Then the motivic sphere $\mathbf{1}[p^{-1}]$ is $f$-pure.
\end{lemma}

\proof
Let $C=i^\triangle(\mathbf{1}[p^{-1}])$. Since $C/\eta=i^\triangle((\mathbf{1}/\eta)[p^{-1}])=0$, it suffices to show that $C[\eta^{-1}]=0$. We know that the $(2,\eta)$-inverted sphere $\mathbf{1}[2^{-1},\eta^{-1}]$ agrees with the minus part of the sphere spectrum $\mathbf{1}_{-}[2^{-1}]$ (\cite[Remark 3.6]{ALP}), which is $i$-pure by \cite[Thm. 3.6]{DFJK}. It follows that $C[2^{-1},\eta^{-1}]=0$. If $p=2$, this proves $C[\eta^{-1}]=0$ since $2$ is already invertible in $C$; if $p$ is odd, the result follows from $C[\eta^{-1}]_{(2)}=i^\triangle(\mathbf{1}[\eta^{-1},p^{-1}])_{(2)}=0$.
\endproof



\subsection{}
By Proposition~\ref{prop:Zpfin} and Lemma~\ref{lm:1pred}, it suffices to prove purity for the spectra $(\mathbf{1}/\eta)[p^{-1}]$, and $\mathbf{1}[\eta^{-1}]_{(2)}$ if $p$ is odd, for closed immersions between regular schemes of finite type over $\mathbb{Z}_{(p)}$. 



\section{Purity for the $\eta$-quotient sphere}

\subsection{}
We first reduce purity for the $\eta$-quotient sphere $(\mathbf{1}/\eta)[p^{-1}]$ to purity for $\mathbf{MGL}[p^{-1}]$:
\begin{proposition}
\label{prop:HF2pur}
Let $p$ be a prime number and let $i:Z\to X$ be a closed immersion between regular schemes of finite dimension with $\operatorname{char}(X)\subset\{0,p\}$. 
If $\mathbf{MGL}[p^{-1}]$ is $i$-pure, then
\begin{enumerate}
\item
The spectrum $(\mathbf{1}/\eta)[p^{-1}]$ is $i$-pure.
\item
If $2\notin\operatorname{char}(X)$, the spectrum $\mathbf{H}\mathbb{F}_2$ is $i$-pure.
\end{enumerate}
\end{proposition}
\proof
We know that there is a factorization of the unit
\begin{align}
\mathbf{1}\to \mathbf{1}/\eta\to \mathbf{MGL}
\end{align}
(see \cite[Theorem 3.8]{Hoy}, \cite[Lemma 3.24]{RSO}), and by \cite[Prop. 3.2.14]{Mon}, the spectrum $(\mathbf{1}/\eta)[p^{-1}]$ is $\eta$-complete and connective. Therefore by \cite[Thm. 1.2.1, 1.2.3 and \S 5.1]{Mon}, the $\mathbf{MGL}$-cobar resolution induces an equivalence
\begin{align}
\operatorname{lim}_{[n]\in\Delta}
\left((\mathbf{1}/\eta)[p^{-1}]\otimes\mathbf{MGL}[p^{-1}]^{\otimes(n+1)}\right)
\simeq
(\mathbf{1}/\eta)[p^{-1}]
\end{align}
for both $X$ and $Z$. Since $(\mathbf{1}/\eta)$ is a cellular spectrum, by Lemma~\ref{lm:celpur}, if $\mathbf{MGL}[p^{-1}]$ is $i$-pure, then each term $(\mathbf{1}/\eta)\otimes\mathbf{MGL}^{\otimes(n+1)}[p^{-1}]$ is also $i$-pure, and the first assertion follows. 

If $p$ is odd, 
since $\mathbf{H}\mathbb{F}_2$ is cellular, by Lemma~\ref{lm:celpur}, purity for $(\mathbf{1}/\eta)[p^{-1}]$ implies purity for $(\mathbf{1}/\eta)[p^{-1}]\otimes\mathbf{H}\mathbb{F}_2$. This implies purity for $\mathbf{H}\mathbb{F}_2$ since we have an isomorphism
\begin{align}
(\mathbf{1}/\eta)[p^{-1}]\otimes\mathbf{H}\mathbb{F}_2\simeq \mathbf{H}\mathbb{F}_2\oplus\mathbf{H}\mathbb{F}_2(1)[2]
\end{align}
which proves the second statement.
\endproof
Note that in \cite[Lemma 5.6 and Cor. 5.10]{Ron}, R\"ondigs used the $\mathbf{kw}$-cobar resolution to prove purity for the $\eta$-inverted sphere $\mathbf{1}[\eta^{-1}]$ over a Dedekind scheme.

\subsection{}
To prove purity for $\mathbf{MGL}[p^{-1}]$, we use Gabber's strategy in \cite[Exp. XVI]{ILO}. 
The following lemma allows to descend purity for $\mathbf{MGL}$ along alterations:
\begin{lemma}
\label{lm:MGLpurdes}
Let $h:Y\to X$ be a projective, dominant, generically finite, generically \'etale morphism of integral regular schemes. 
Let $g:X\to S$ be another morphism such that both $g$ and $gh$ are smoothable lci. If $\mathbf{MGL}_S[e^{-1},\operatorname{char}(X)^{-1}]$ is $gh$-pure, then it is also $g$-pure.
\end{lemma}
\proof
Under our assumptions, the morphism $h$ is smoothable lci of relative dimension $0$. Therefore the purity transformation associated to $h$ on $\mathbf{MGL}$ modules gives a natural transformation $h^*\to h^!$ (see \cite{Deg11}, \cite[4.4.3]{DJK}). Denote 
\begin{align}
a_h:\mathbf{MGL}_X\to h_*h^*\mathbf{MGL}_X\to h_*h^!\mathbf{MGL}_X\to \mathbf{MGL}_X
\end{align}
which we consider as a class $a_h\in \mathbf{MGL}^{0,0}(X)$. Let $r$ be the virtual relative dimension of $g$, and denote by $\eta_g:\mathbf{MGL}_X(r)[2r]\to g^!\mathbf{MGL}_S$ the $\mathbf{MGL}$-valued fundamental class. Let $\omega$ be the composition
\begin{align}
\begin{split}
\omega:
g^!\mathbf{MGL}_S
&\to 
h_*h^*g^!\mathbf{MGL}_S
\to 
h_*h^!g^!\mathbf{MGL}_S
\overset{(\eta_{gh})^{-1}}{\simeq}
h_*\mathbf{MGL}_Y(r)[2r]\\
&\to
h_*h^!\mathbf{MGL}_X(r)[2r]
\to
\mathbf{MGL}_X(r)[2r].
\end{split}
\end{align}
Standard unit-counit computations give identities $\omega\circ\eta_g=a_h\cdot Id$ and $\eta_g\circ\omega=a_h\cdot Id$. Therefore to show that $\eta_g$ is an isomorphism, it suffices to show that the class $a_h\in \mathbf{MGL}^{0,0}(X)[e^{-1},\operatorname{char}(X)^{-1}]$ is invertible. 

By Lemma~\ref{lm:ptcons}, it suffices to show that the restriction of $a_h$ to every point $x$ of $X$ is invertible in $\mathbf{MGL}^{0,0}(\kappa(x))[e^{-1}]$. Let $p$ be the exponential residual characteristic of $x$. By \cite[Prop. 8.2]{Hoy}, we have an isomorphism $\mathbf{MGL}^{0,0}(\kappa(x))[p^{-1}]\simeq\mathbb{Z}[1/p]$. We claim that for every point $x$, the restriction of $a_h$ at $x$ is equal to $e$. Indeed, the orientation of $\mathbf{KGL}$ induces a map $\mathbf{MGL}\to\mathbf{KGL}$, and consequently we have maps
\begin{align}
\label{eq:MGLKGLX}
\mathbf{MGL}^{0,0}(X)[\operatorname{char}(X)^{-1}]\to \mathbf{KGL}^{0,0}(X)[\operatorname{char}(X)^{-1}],
\end{align}
\begin{align}
\label{eq:MGLKGLpt}
\mathbf{MGL}^{0,0}(\kappa(x))[p^{-1}]\simeq \mathbf{KGL}^{0,0}(\kappa(x))[p^{-1}]\simeq\mathbb{Z}[1/p]
\end{align}
where the second isomorphism of~\eqref{eq:MGLKGLpt} is induced by the rank map.
If $f$ is finite \'etale over $x$, a trace computation shows that the restriction of $a_h$ at $x$ has rank $e$, see \cite[Exp. XVI, Lemme 3.5.10]{ILO}. In general, since $f$ is generically finite \'etale, such $x$ form a dense open subset of $X$. On the other hand, the image of $a_h$ under the map~\eqref{eq:MGLKGLX} has locally constant (virtual) rank. We conclude that the restriction of $a_h$ at every point $x$ has rank $e$, which finishes the proof.
\endproof

Note that Lemma~\ref{lm:MGLpurdes} also applies to $\mathbf{H}\mathbb{Z}$ in the place of $\mathbf{MGL}$, but not to the sphere spectrum. This is the only place in the proof, apart from cellularity, where a specific property of the spectrum $\mathbf{MGL}$ is used. 

\subsection{}
\label{num:DVRpur}
Recall that if $i$ is the inclusion of a closed point in a DVR of mixed characteristic $(0,p)$, the sphere spectrum $\mathbf{1}[p^{-1}]$ is $i$-pure (\cite[Thm. 7.8]{RS}). Since $\mathbf{MGL}$ is cellular, we see that $\mathbf{MGL}[p^{-1}]$ is also $i$-pure.




\subsection{}
In the following lemma we study torus actions on affine spaces and some properties of equivariant objects.
\begin{lemma}
\label{lm:eqret}
\begin{enumerate}
\item
Let $B$ be a scheme and denote by $q:\mathbb{A}^r_B\to B$ the projection with zero section $s$.
Consider a $\mathbb{G}_m$-action on $\mathbb{A}^r_B$ of weights $a_1,\cdots,a_r$ such that all the $a_i$'s are nonzero and have the same sign. If $C\in\mathbf{SH}(\mathbb{A}^r_B)$ is an equivariant object, then the adjunction induces an isomorphism
\begin{align}
\label{eq:qscont}
q_*C\simeq s^*C.
\end{align}
\item
Consider the affine plane $\mathbb{A}^2_B$ with coordinates $x,y$. Let $m,n>0$ be integers and consider the $\mathbb{G}_m$-action on $\mathbb{A}^2_B$ given by $\lambda\cdot(x,y)=(\lambda^mx,\lambda^{-n}y)$. We have a Cartesian diagram
\begin{align}
\begin{split}
  \xymatrix@=10pt{
Z\bigcap U=\{x\ne0,y=0\} \ar[r]^-{j'} \ar[d]_-{i'} & Z=\{y=0\} \ar[d]^-{i} \\
U=\{x\ne0\} \ar[r]^-{j} & \mathbb{A}^2_B. 
  }
\end{split}
\end{align}
Let $N\in\mathbf{SH}(U)$ be an equivariant object.
\begin{enumerate}
\item
If $i'^*N=0$, then $i^*j_*N=0$.
\item
In general, the adjunction induces an isomorphism 
\begin{align}
\label{eq:ijex}
i^*j_*N\simeq j'_*i'^*N.
\end{align}
\end{enumerate}
\end{enumerate}
\end{lemma}
\proof
\begin{enumerate}
\item
Replacing $C$ by the fiber of $C\to s_*s^*C$, we are reduced to show that if $s^*C=0$, then $q_*C=0$. We may assume that all the $a_i$'s are nonnegative, so that the action extends to a contraction $\phi:\mathbb{A}^1\times\mathbb{A}^r_B\to\mathbb{A}^r_B$. Denote by $p:\mathbb{A}^1\times\mathbb{A}^r_B\to\mathbb{A}^r_B$ the projection, and by equivariance we have $(\phi^*C)_{|\{0\}\times\mathbb{A}^r_B}=0$ and $(\phi^*C)_{|\mathbb{G}_m\times\mathbb{A}^r_B}=(p^*C)_{|\mathbb{G}_m\times\mathbb{A}^r_B}$. By localization and $\mathbb{A}^1$-invariance we obtain $p_*\phi^*C=0$, and therefore $q_*\phi_*\phi^*C=q_*p_*\phi^*C=0$.
The unit section of $\phi$ exhibits $q_*C$ as a retract of $q_*\phi_*\phi^*C$, which finishes the proof.
\item
The second claim can be obtained by applying the first claim to the fiber of$N\to i'_*i'^*N$, so it suffices to prove the first one. Consider the compactification on the $x$-coordinate $\pi:P=\mathbb{P}^1_x\times_B\mathbb{A}^1_y\to B$ of $\mathbb{A}^2_B$, with an open chart $l:W=\mathbb{A}^1_z\times_B\mathbb{A}^1_y\to P$ where $z=x^{-1}$. Let $\widetilde{N}$ be the extension by zero of $N$ to $W$, and let $C=l_*\widetilde{N}\in\mathbf{SH}(P)$. The $\mathbb{G}_m$-action on $\mathbb{A}^2_B$ extends to $P$, and on the chart $W$ the coordinates $z$ and $y$ have weights $-m$ and $-n$. By~\eqref{eq:qscont} we have $\pi_*C\simeq\widetilde{N}_{|\{z=y=0\}}=0$. We have a Cartesian diagram
\begin{align}
\begin{split}
  \xymatrix@=10pt{
P \ar[r]^-{v} \ar[d]_-{} & \mathbb{A}^1_y \ar[d]^-{d} \\
\mathbb{P}^1_x \ar[r]^-{w} & B
  }
\end{split}
\end{align}
and using~\eqref{eq:qscont} again we obtain $0=\pi_*C=d_*v_*C\simeq(v_*C)_{|\{y=0\}}\simeq w_*C_{|\{y=0\}}$. 

On the other hand, since $i'^*N=0$ and $\widetilde{N}$ is extended to $W$ by zero at infinity, we know that the object $C_{|\{y=0\}}\in\mathbf{SH}(\mathbb{P}^1_x)$ vanishes on $W\bigcap\{y=0\}$. It follows that $C_{|\{y=0\}}$ is supported on $\{x=0\}$, and therefore $C_{|\{y=0\}}=0$. We conclude as $i^*j_*N=(C_{|\{y=0\}})_{|Z}=0$.

\end{enumerate}
\endproof

\subsection{}
Let $D$ be an excellent DVR of mixed-characteristic $(0,p)$, and fix a uniformizer $\pi$ of $D$. Let $S=\operatorname{Spec} D$ with generic point $\eta=\operatorname{Spec}K$ and closed point  $s=\operatorname{Spec}k$. 
\begin{lemma}
\label{lm:chartpur}
Let $e_1,\cdots,e_r\geqslant1$ be integers and consider the morphism 
\begin{align}
f:X=\operatorname{Spec} D[t_1,\cdots,t_r]/(t_1^{e_1}\cdots t_r^{e_r}-\pi)\to S. 
\end{align}
For $a=1,\cdots,r$, consider the morphism $i_a:H_a=V(t_a)\to X$. By \cite[Exp. XVI, Lemme 3.5.13 (iii)]{ILO}, $X$ is regular and $H_a\simeq\mathbb{A}^1_k$ is a regular Cartier divisor.
\begin{enumerate}
\item
\label{num:chartpur1}
Write $J:X_K\to X$ the inclusion. For every subset $I\subset\{1,\cdots,r\}$, denote $H_I=\bigcup_{a\in I}H_a$, $H_I^\circ=H_I\setminus\bigcup_{a\notin I}H_a$ and $h_I:H_I\to X$, $v_I:H_I^\circ\to H_I$ the inclusions. Let $M\in \mathbf{SH}(K)$.
Then 
the adjunction induces an isomorphism 
\begin{align}
\label{eq:adjiso}
h_I^*J_*M_{|X_K}\simeq v_{I,*}v_I^*h_I^*J_*M_{|X_K}.
\end{align}
\item
\label{num:chartpur2}
The sphere $\mathbf{1}[p^{-1}]$ is $f$-pure and $i_a$-pure for $a=1,\cdots,r$.
\end{enumerate}
\end{lemma}
\proof
\begin{enumerate}
\item
If $I=\emptyset$ or $I=\{1,\cdots,r\}$, the inclusion $v_I$ is an isomorphism and the result is clear. Fix $a\notin I$ and denote $j_a:D(t_a)\to X$, $k_a:H_I\bigcap D(t_a)\to D(t_a)$. We will show that we can remove the components $H_a$ from $H_I$ one by one. 

For $b\in I$, regard $X$ as a closed subscheme of the plane $\mathbb{A}^2\times \operatorname{Spec} D[t_i]_{i\ne a,b}$, with coordinates $x=t_a$ and $y=t_b$. Consider the $\mathbb{G}_m$-action where $t_a$ has weight $e_b$ and $t_b$ has weight $-e_a$. Let $N\in \mathbf{SH}(D(t_a))$ be equivariant for each of these actions indexed by $b\in I$, and let $Q$ be the fiber of $N\to k_{a,*}k_a^*N$. Then $Q$ restricts to zero on every $H_b\bigcap D(t_a)$. By~\eqref{eq:ijex}, $j_{a,*}Q$ restricts to zero on every $H_b$, and therefore by localization we have $h_I^*j_{a,*}Q=0$. This shows that the isomorphisms persists after removing $H_a$ from $H_I$. Repeating the process for every $a\notin I$ gives the result.



\item
We first show that $\mathbf{1}[p^{-1}]$ is $i_a$-pure.
If $r=1$, this follows from~\ref{num:DVRpur}. Now assume $r\geqslant2$ and fix $a\in\{1,\cdots,r\}$. Denote $Y=\bigcup_{b\ne a}H_a$ and $U=X\setminus H_a$, where $Y$ is endowed with the reduced scheme structure. We have a commutative diagram where both squares are Cartesian and $\operatorname{Tor}$-independent:
\begin{align}
\begin{split}
  \xymatrix@=10pt{
   & Y\bigcap U \ar[r]^-{v} \ar[d]_-{d} & Y \ar[d]^-{h} & H_a\bigcap Y \ar[l]_-{c} \ar[d]^-{s} \\
   X_K \ar[r]^-{u} &U \ar[r]^-{j} & X & H_a. \ar[l]_-{i_a}
  }
\end{split}
\end{align}
By~\eqref{eq:adjiso} the adjunction induces an isomorphism $h^*j_*u_*\mathbf{1}[p^{-1}]\simeq v_*d^*u_*\mathbf{1}[p^{-1}]$. Since $X_K\bigcap Y=\emptyset$, the cofiber of $\mathbf{1}[p^{-1}]\to u_*\mathbf{1}[p^{-1}]$ has the form $d_*\phi$, from which we deduce an isomorphism $h^*j_*\mathbf{1}[p^{-1}]\simeq v_*\mathbf{1}[p^{-1}]$. By localization triangle we have an isomorphism $s^*i_a^!\mathbf{1}[p^{-1}]\simeq c^!h^*\mathbf{1}[p^{-1}]$.

Let $C_a=i_a^{\triangle}(\mathbf{1}_X[p^{-1}])\in \mathbf{SH}(H_a)$. Note that $Y\simeq (H_a\bigcap Y)\times\mathbb{A}^1$ with the morphism $c$ identified as the zero section, so $s^*C_a=0$.
Let $O_a=H_a\setminus Y=H_a\setminus\bigcup_{b\ne a}H_a$, and we are reduced to show that $C_{a|O_a}=0$.

Let $\xi_a$ be the generic point of $O_a$. The local ring $\mathcal{O}_{X,\xi_a}$ is a DVR of mixed characteristic $(0,p)$, so by~\ref{num:DVRpur} we have $(C_a)_{|\xi_a}=0$. By Lemma~\ref{lm:ptcons}, it remains to show that for any point $z$ of $O_a$ we have $C_a|_z=0$.
The torus $G=(\Gm^{r-1})_D$ with coordinates $(u_j)_{j\ne a}$ acts on $X$ by $t_j\mapsto u_j^{e_a}t_j$ for $j\ne a$ and $t_a\mapsto(\prod_{j\ne a}u_j^{-e_j})t_a$. The action preserves $H_a$, and both the restricted action and projection maps $\mu_a,p_a:G\times_DH_a\to H_a$ are smooth, and we have $\mu^*_aC_a\simeq p^*_aC_a$. For $z\in O_a$, translate $z$ by the generic element of $G_{\kappa(z)}$: over the purely transcendental field extension $\kappa(z)(u_j)_{j\ne a}$ of $\kappa(z)$, its coordinates $t_j(z)u_j^{e_a}$ are algebraically independent, so the translated point factors through $\xi_a$. Lemma~\ref{lm:extcons} shows that $C_a|_z=0$, and it follows that $\mathbf{1}[p^{-1}]$ is $i_a$-pure.

Now let $C=f^{\triangle}(\mathbf{1}_S[p^{-1}])$. Since the generic fiber of $f$ is smooth, we have $C_{|X_K}=0$. For each $a$, the composite $fi_a$ factors as $H_a\to\operatorname{Spec}k\to S$, so by~\ref{num:DVRpur}, $\mathbf{1}[p^{-1}]$ is $fi_a$-pure. This shows that $i_a^!C=0$. Denote by $\phi_1:X\setminus H_1\to X$ the inclusion, then $C\simeq \phi_{1,*}\phi_1^*C$; repeating the argument with $H_2,\cdots,H_r$, we are reduced to $C_{|X_K}$ which is $0$. Therefore we have proved that $\mathbf{1}[p^{-1}]$ is $f$-pure.
\end{enumerate}
\endproof
\begin{remark}

Part~\eqref{num:chartpur1} of the lemma is analogous to \cite[Lemma 3.1]{Fal}, \cite[Lemma 2.2.1]{KS} and \cite[Thm. 3.3.10 and 3.3.24]{Ayo}, and while part~\eqref{num:chartpur2} generalizes \cite[Exp. XVI, Lemme 3.5.13 (iv)]{ILO}. Part~\eqref{num:chartpur1} can also be proven using \cite[Lemma 2.23]{Pha}.
\end{remark}


\subsection{}
The following lemma is an analogue of \cite[Exp. XVI, Th. 3.6.2]{ILO}, where we avoid using \'etale neighborhoods:
\begin{lemma}
\label{lm:Gabber362}
Let $g:Y\to S$ be a quasi-projective morphism such that
\begin{enumerate}
\item
The source $Y$ is regular.
\item
The reduced special fiber $(Y_s)_{\operatorname{red}}$ is a simple normal-crossing divisor.
\item
The induced morphism $(Y,(Y_s)_{\operatorname{red}})\to (S,s)$ is log smooth.
\end{enumerate}
Then the sphere $\mathbf{1}_S[p^{-1}]$ is $g$-pure. Consequently $\mathbf{MGL}_S[p^{-1}]$ is also $g$-pure.
\end{lemma}
\proof
The generic fiber of $g$ is smooth. Let $y\in Y_s$ be a point in the special fiber, and let $U\subset Y$ be a Zariski open neighborhood of $y$. Shringking $U$, we may assume that any irreducible component of $(Y_s)_{\operatorname{red}}$ that meets $U$ passes through $y$. Let $Z_1,\cdots,Z_r$ be these components, with multiplicity $e_1,\cdots,e_r$. Choose a local equation $t_i$ on $U$ defining each $Z_i$. Then we have $(g_{|U})^*\pi=ut_1^{e_1}\cdots t_r^{e_r}$ for some unit $u\in\mathcal{O}^\times(U)$ (\cite[Prop. 2.5]{AIS}). It suffices to show that $\mathbf{1}_S[p^{-1}]$ is pure relatively to the morphism $g:U\to S$.

Define a $D$-scheme
\begin{align}
W=\operatorname{Spec} D[x_1,\ldots,x_r,v^{\pm1}]/(vx_1^{e_1}\cdots x_r^{e_r}-\pi)
\end{align}
and a morphism 
\begin{align}
f:U\times\mathbb{G}_m^r
=
\operatorname{Spec}_{\mathcal{O}_U}\mathcal{O}_U[z_1^{\pm1},\cdots,z_r^{\pm1}]
\to W
\end{align}
sending $x_i$ to $z_it_i$ and $v$ to $uz_1^{-e_1}\cdots z_r^{-e_r}$. We claim that $f$ is smooth. Indeed, endow each of the schemes $S,Y,W,U,U\times\mathbb{G}_m^r$ with the log structure arising from the special fiber.
The chart on $W$ has generators $x_1,\cdots,x_r$ with log parameter $vx_1^{e_1}\cdots x_r^{e_r}$; pulling them back along the morphism $U\to W$ given by $(t_1,\cdots, t_r,u)$ gives the generators $t_1,\cdots,t_r$ of chart on $U$ with log parameter $ut_1^{e_1}\cdots t_r^{e_r}$ (\cite[3.1.2]{ATW}). Therefore the induced map $U \to [W/\mathbb{G}_m^r]\to\mathcal{L}og_S$ 
classfies the log morphism associated to $g_{|U}:U\to S$ (\cite{Ols}), and we have a Cartesian diagram
\begin{align}
\begin{split}
  \xymatrix@=10pt{
    U\times\mathbb{G}_m^r\ \ar[r]^-{f} \ar[d]_-{} & W \ar[d]^-{} \\
    U \ar[r]_-{} & [W/\mathbb{G}_m^r].
  }
\end{split}
\end{align} 
Since $g:Y\to S$ is log smooth, by \cite[Thm. 4.6, Cor. 5.31]{Ols} we know that $U \to [W/\mathbb{G}_m^r]$ is smooth, which implies the smoothness of $f$.

Now $W$ is an open subscheme of $\operatorname{Spec} D[x_1,\ldots,x_r,v]/(vx_1^{e_1}\cdots x_r^{e_r}-\pi)$. By Lemma~\ref{lm:chartpur}, we know that $\mathbf{1}_S[p^{-1}]$ is pure relatively to the morphism $W\to S$. It follows that $\mathbf{1}_S[p^{-1}]$ is pure relatively to the morphism $U\times\mathbb{G}_m^r\to S$. Since the projection $U\times\mathbb{G}_m^r\to U$ is smooth with a section, this implies that $\mathbf{1}_S[p^{-1}]$ is pure relatively to the morphism $U\to S$.
\endproof


\begin{proposition}
\label{prop:dvrpur}
Let $q:X\to S$ be a quasi-projective morphism with $X$ regular. Then $\mathbf{MGL}[p^{-1}]$ is $q$-pure.
\end{proposition}
\proof
Any quasi-projective morphism between regular schemes is smoothable lci. We may also assume that $q$ is surjective. By \cite[Thm. 4.3.1]{Tem}, there exists a commutative square
\begin{align}
\begin{split}
  \xymatrix@=10pt{
    Y \ar[r]^-{h} \ar[d]_-{g} & X \ar[d]^-{q} \\
    S' \ar[r]^-{f} & S.
  }
\end{split}
\end{align} 
where $Y$ and $S'$ are regular, the horizontal maps are projective $p$-alterations, the inverse image of the special fiber $h^{-1}(X_s)$ is contained in a simple normal-crossing divisor, and $Y\to S'$ is log smooth and quasi-projective. In addition, $S'$ and $S$ are integral curves, so the morphism $f$ is finite and generically \'etale. 

We claim that $\mathbf{MGL}[p^{-1}]$ is $f$-pure. To see this, let $s'$ be a closed point of $S'$ lying over a point $i:s\to S$. By~\ref{num:DVRpur}, $\mathbf{MGL}[p^{-1}]$ is $i$-pure, and therefore $\mathbf{MGL}[p^{-1}]$ is pure relatively to $s'\to S$. By \cite[Prop. B.20]{DFJK} this implies that $\mathbf{MGL}[p^{-1}]$ is $f$-pure. Now by Lemma~\ref{lm:Gabber362}, we know that $\mathbf{MGL}[p^{-1}]$ is $g$-pure, so it is also $qh=fg$-pure; by Lemma~\ref{lm:MGLpurdes} we conclude that $\mathbf{MGL}[p^{-1}]$ is $q$-pure.
\endproof

\begin{corollary}
\label{cor:MGLpur}
Let $D$ be an excellent DVR of mixed-characteristic $(0,p)$, and let $i:Z\to X$ be a closed immersion between regular schemes of finite type over $D$. Then $\mathbf{MGL}[p^{-1}]$ is $i$-pure.
\end{corollary}
\proof
Purity can be checked Zariski locally, and therefore we may assume both $X$ and $Z$ affine, which are therefore quasi-projective over $D$. Then the result follows from Lemma~\ref{lm:DJK4310} and Proposition~\ref{prop:dvrpur}.
\endproof

\subsection{}
We end this section with a discussion on the purity for $\mathbf{MGL}$ at the characteristic, using the recent work of Bachmann-Hopkins (\cite{BH}), which is not used in the proof of Theorem~\ref{thm:pur1}.
\begin{lemma}
\label{lm:purMGLHZ}
Let $f$ be a smoothable lci morphism.
Then $\mathbf{MGL}$ is $f$-pure if and only it $\mathbf{H}\mathbb{Z}$ is $f$-pure.
\end{lemma}
Indeed, the integral Hopkins-Morel isomorphism (\cite[Thm. 7.2]{BH}) induces an isomorphism $\mathbf{MGL}/(t_1,t_2,\cdots)\simeq\mathbf{H}\mathbb{Z}$. Since purity is preserved by colimits (\cite[Lemma 7.9]{RS}), purity for $\mathbf{MGL}$ implies purity for $\mathbf{H}\mathbb{Z}$. Conversely, by \cite[Thm. 7.3]{BH} the slice filtration on $\mathbf{MGL}$ has $n$-th graded piece $s_n\mathbf{MGL}\simeq\mathbf{H}\mathbb{Z}(n)[2n]\otimes L_n$, where $L_n$ is the degree $n$ part of the Lazard ring. Using convergence of the slice tower (obtained by combining \cite[Lemma 2.42]{BBX} and \cite[Thm. 7.3]{BH}), we see that purity for $\mathbf{H}\mathbb{Z}$ implies purity for $\mathbf{MGL}$.

As a consequence, since purity of $\mathbf{H}\mathbb{Z}$ is known over Dedekind domains (\cite[Thm. 2.13]{AE}, \cite[Thm. 7.6]{RS}, \cite[Cor. 8.4]{FS}), we deduce purity for $\mathbf{MGL}$ in that case:
\begin{corollary}
\label{cor:MGLDpur}
Let $i:Z\to X$ be the inclusion of a closed point into the spectrum of a Dedekind domain. Then $\mathbf{MGL}$ is $i$-pure.
\end{corollary}

\section{Purity for the $2$-local $\eta$-inverted sphere}

\subsection{}
\label{num:Bacsp}
In this section, we work over $\mathbb{Z}[1/2]$-schemes. Consider the following motivic spectra in $\mathbf{SH}(\operatorname{Spec}\mathbb{Z}[1/2])$:
\begin{enumerate}
\item
$\mathbf{K}^W$ the spectrum representing powers of the fundamental ideal $I^*$ (\cite[Def. 2.6]{Bac22}),
\item
$\mathbf{KQ}$ the Hermitian K-theory spectrum, and $\mathbf{KW}=\mathbf{KQ}[\eta^{-1}]$ the periodic Witt spectrum, with connective cover $\mathbf{kw}$ (\cite[Def. 1.3]{ALP}, \cite[\S1]{Bac22}),
\item
$\mathbf{HW}$ the Witt cohomology spectrum (\cite[5.8]{DFJK}, \cite[Lemma 3.9]{Bac22}).
\end{enumerate}
For $E$ any of these spectra and for any $\mathbb{Z}[1/2]$-scheme $S$, denote by $E_S\in\mathbf{SH}(S)$ the restriction to $S$ (\cite[Def. 4.6]{Bac22}). These spectra are endowed with the following isomorphisms:
\begin{align}
\label{eq:KWrel}
\mathbf{K}^W/\eta\simeq\mathbf{H}\mathbb{F}_2/\tau,
\qquad
\mathbf{K}^W[\eta^{-1}]\simeq\mathbf{HW},
\qquad
\mathbf{kw}/\beta\simeq\mathbf{HW},
\qquad
\mathbf{kw}[\beta^{-1}]\simeq\mathbf{KW},
\end{align}
where $\tau\in\mathbf{H}\mathbb{F}_2^{0,1}(\operatorname{Spec}\mathbb{Z}[1/2])$ is the unique nonzero element (\cite[Example 2.4]{Bac22}), and $\beta\in\mathbf{kw}^{-4,0}(\operatorname{Spec}\mathbb{Z}[1/2])$ is the Bott class for connective K-theory (\cite{Ana}, where the class denoted as $\beta$ is $\eta^4\beta$ with our notations). Indeed, the first three isomorphisms of~\eqref{eq:KWrel} are \cite[Lemmas 2.7, 3.9 and 4.4]{Bac22}, and the last one follows from the computation of homotopy sheaves in \cite[Remark 3.7 and Lemma 4.4]{Bac22}.


\subsection{}
Let $T$ be a noetherian scheme. 
By \cite[Thm. 2.2]{Voe}, every cdh cover of $T$ refines to a finite iteration of elementary Nisnevich squares and abstract blowup squares of the form
\begin{align}
\label{eq:elemsq}
\begin{split}
  \xymatrix@=10pt{
    E \ar[r]^-{} \ar[d]_-{} & Y \ar[d]^-{} \\
    Z \ar[r]_-{} & T.
  }
\end{split}
\end{align} 

\begin{lemma}
\label{lm:cdhdim}
If $T$ has dimension $d$, one may choose these squares such that the resulting cdh cover $(U_i)_{i\in I}$ of $T$ has each member $U_i$ of dimension $\leqslant d$.
\end{lemma}
\proof
It suffices to show that these squares can be chosen such that all four schemes in each square has dimension at most $d$. For Nisnevich squares this follows from \'etaleness. Consider an abstract blowup square of the form~\eqref{eq:elemsq}. Let $Y'$ be the scheme-theoretic closure of $Y\times_T(T-Z)$ in $Y$. Then $Y'$ is still proper over $T$ and induces an isomorphism $Y'\times_T(T-Z)\simeq T-Z$. The irreducible components of $Y'$ are birational to the irreducible components of $T_{\operatorname{red}}$, so 
they have dimension at most $d$. We replace $Y$ by $Y'$, and by induction we apply the same replacement for every square obtained subsequently by pullbacks, which finishes the proof.
\endproof

\subsection{}
Let $p$ be an odd prime and let $D=\mathbb{Z}_{(p)}[u]/(u^2+p-1)$. Then $D$ is finite \'etale over $\mathbb{Z}_{(p)}$ since the discriminant $4(1-p)$ is a unit; its fraction field is the imaginary quadratic field $\mathbb{Q}[\sqrt{1-p}]$, and its special fiber is $\mathbb{F}_p\times\mathbb{F}_p$. Consequently the family $\{\operatorname{Spec}\mathbb{Q},\operatorname{Spec} D\}$ form a Nisnevich cover of $\operatorname{Spec}\mathbb{Z}_{(p)}$. We know that no residue field of $D$ is orderable, and therefore the same holds for any $D$-scheme $V$, and the associated real scheme $V_r$ is empty. 
\begin{lemma}
\label{lm:Ijfin}
Let $V$ be a $D$-scheme of finite type. Then for any $m>\operatorname{dim}(V)+2$, we have
\begin{align}
R\Gamma_{\operatorname{cdh}}(V,\underline{I}_{\operatorname{cdh}}^m)=0.
\end{align}
Here $\underline{I}_{\operatorname{cdh}}^m$ is the cdh sheaf of the $m$-th power of the fundamental ideal, with $\underline{I}^m_{\operatorname{cdh}}=\underline{W}_{\operatorname{cdh}}$ for $m\leqslant0$.
\end{lemma}
\proof
Let $d=\operatorname{dim}(V)$, and let $T$ be a $D$-scheme of finite type of dimension at most $d$. By Gabber rigidity of Witt rings with respect to henselian local rings (\cite[Lemma 5.2]{Jac2}, \cite[Prop. 10.5 and 10.19]{BEM}), the stalk of the Nisnevich sheaf $\underline{I}_{\operatorname{Nis}}^m$ at a point $x$ is $I^m(\kappa(x))$. The residue field $\kappa(x)$ is not orderable and has transcendence degree at most $d$ over the corresponding residue field of $D$, which is either a finite field or a imaginary quadratic extension of $\mathbb{Q}$. By \cite[\S2]{Ser}, we have $\operatorname{cd}_2\kappa(x)\leqslant d+2$. So by \cite[Prop. 2.3]{Bac22}, for $m>d+2$ we have $I^m(\kappa(x))=0$, and therefore $\underline{I}_{\operatorname{Nis}}^m(T)=0$.

It remains to pass from $\underline{I}_{\operatorname{Nis}}^m$ to $\underline{I}_{\operatorname{cdh}}^m$. We know that cdh locally, a section $\sigma$ of $\underline{I}_{\operatorname{cdh}}^m(T)$ lifts to a section of $\underline{I}_{\operatorname{Nis}}^m(T)$. By Lemma~\ref{lm:cdhdim}, we can find a finite cdh cover $(U_i)_{i\in I}$ of $T$ with $\operatorname{dim}(U_i)\leqslant d$ such that $\sigma_{U_i}=0$. This proves that $\underline{I}_{\operatorname{cdh}}^m(T)=0$.

We now show by induction that for every $q\geqslant0$ and every $T$ separated of finite type over $D$ of dimension $\leqslant d$ we have $H^q_{\operatorname{cdh}}(T,\underline{I}_{\operatorname{cdh}}^m)=0$. The case $q=0$ is already proven. Assume that $q\geqslant1$, and the statement holds for $q-1$ and every such $T$. Let $\sigma\in H^{q}_{\operatorname{cdh}}(T,\underline{I}_{\operatorname{cdh}}^m)$. We know that there exists a cdh cover $(T_i)$ of $T$ such that $\sigma_{|T_i}=0$. By Lemma~\ref{lm:cdhdim}, we can find a finite cdh cover $(U_i)_{i\in I}$ of $T$ with $\operatorname{dim}(U_i)\leqslant d$ such that $\sigma_{U_i}=0$. On the other hand, for any cdh square of the form~\eqref{eq:elemsq}, by \cite[Lemma 2.19]{Voecd} we have an exact sequence
\begin{align}
H^{q-1}_{\operatorname{cdh}}(E,\underline{I}_{\operatorname{cdh}}^m)
\to
H^{q}_{\operatorname{cdh}}(T,\underline{I}_{\operatorname{cdh}}^m)
\to
H^{q}_{\operatorname{cdh}}(Z,\underline{I}_{\operatorname{cdh}}^m)
\oplus
H^{q}_{\operatorname{cdh}}(Y,\underline{I}_{\operatorname{cdh}}^m).
\end{align}
By induction hypothesis we have 
$H^{q-1}_{\operatorname{cdh}}(E,\underline{I}_{\operatorname{cdh}}^m)=0$, 
so $H^{q}_{\operatorname{cdh}}(T,\underline{I}_{\operatorname{cdh}}^m)$ injects into 
$H^{q}_{\operatorname{cdh}}(Z,\underline{I}_{\operatorname{cdh}}^m)
\oplus
H^{q}_{\operatorname{cdh}}(Y,\underline{I}_{\operatorname{cdh}}^m)$. Repeating the argument iteratedly with the cdh squares in Lemma~\ref{lm:cdhdim}, it follows that $H^{q}_{\operatorname{cdh}}(T,\underline{I}_{\operatorname{cdh}}^m)$ injects into $\oplus_iH^{q}_{\operatorname{cdh}}(U_i,\underline{I}_{\operatorname{cdh}}^m)$. It follows that $\sigma=0$, which finishes the induction step.

It follows that for every $T$ separated of finite type over $D$ of dimension $\leqslant d$, we have $R\Gamma_{\operatorname{cdh}}(T,\underline{I}_{\operatorname{cdh}}^m)=0$. A \v Cech descent argument shows that the same statement also holds for all $T$ of finite type over $D$ of dimension $\leqslant d$, including the case $T=V$, which finishes the proof. 
\endproof




\begin{remark}
In general, vanishing of Nisnevich cohomology does not imply vanishing of cdh cohomology, see \cite[Exercise 12.32]{MVW}.
\end{remark}



\begin{lemma}
\label{lm:KWcomp}
If $U$ is a $D$-scheme of finite type, then $\mathbf{K}^W_U$ is $\eta$-complete.
\end{lemma}
\proof
Let $V$ be a smooth $U$-scheme. By \cite[Prop. 10.19]{BEM}, for any integers $n$ and $j$ we have an equivalence
\begin{align}
\operatorname{Map}_{\mathbf{SH}(U)}(\Sigma_U^\infty V_+,\mathbf{K}^W_U(n+j)[n+2j])
\simeq
R\Gamma_{\operatorname{cdh}}(V,\underline{I}_{\operatorname{cdh}}^{j+n})[j]
\end{align}
which vanishes if $n+j>\operatorname{dim}(V)+2$ by Lemma~\ref{lm:Ijfin}. Thus the inverse
limit of the tower $\mathbf{K}^W_U(n)[n]$ with transition maps induced by $\eta$, is zero. Taking limits over the cofibers gives an isomorphism $\mathbf{K}^W_U\simeq\operatorname{lim}_n\mathbf{K}^W_U/\eta^n$.
\endproof

\begin{proposition}
\label{prop:1eta2pur}
Let $p$ be an odd prime. Let $i:Z\to X$ be a closed immersion between regular schemes of finite type over $\mathbb{Z}_{(p)}$. Then the $2$-local $\eta$-inverted sphere $\mathbf{1}[\eta^{-1}]_{(2)}$ is $i$-pure.
\end{proposition}
\proof
By Lemma~\ref{lm:purbc}, it suffices to show that $\mathbf{1}[\eta^{-1}]_{(2)}$ is pure relatively to the two morphisms
\begin{align}
i_{\mathbb{Q}}:
Z\times_{\operatorname{Spec}\mathbb{Z}_{(p)}}\operatorname{Spec}\mathbb{Q}
\to
X\times_{\operatorname{Spec}\mathbb{Z}_{(p)}}\operatorname{Spec}\mathbb{Q},
\end{align}
\begin{align}
i_{D}:
Z\times_{\operatorname{Spec}\mathbb{Z}_{(p)}}\operatorname{Spec}D
\to
X\times_{\operatorname{Spec}\mathbb{Z}_{(p)}}\operatorname{Spec}D.
\end{align}
Since $i_{\mathbb{Q}}$ is a closed immersion between regular schemes of equal characteristic $0$, its purity follows from \cite[Thm. C.1]{DFJK}. 

It remains to prove purity for $i_D$. We will be using the four isomorphims in~\eqref{eq:KWrel}. By Proposition~\ref{prop:HF2pur} and Corollary~\ref{cor:MGLpur} we know that $\mathbf{H}\mathbb{F}_2$ is $i_D$-pure, and from the isomorphism $\mathbf{K}^W/\eta\simeq\mathbf{H}\mathbb{F}_2/\tau$, we deduce that $\mathbf{K}^W/\eta$ is $i_D$-pure. By taking successive $\eta$-quotients, we know that each $\mathbf{K}^W/\eta^n$ is $i_D$-pure. By Lemma~\ref{lm:KWcomp}, we deduce that $\mathbf{K}^W$ is $i_D$-pure. From the isomorphisms $\mathbf{K}^W[\eta^{-1}]\simeq\mathbf{HW}\simeq\mathbf{kw}/\beta$ , we deduce that $\mathbf{kw}/\beta$ is $i_D$-pure. 
On the other hand, by \cite[Thm. D.1]{DFJK}, $\mathbf{KQ}$ is $i_D$-pure, and therefore so is $\mathbf{kw}[\beta^{-1}]\simeq\mathbf{KW}=\mathbf{KQ}[\eta^{-1}]$. We deduce that $\mathbf{kw}$ is $i_D$-pure, and from the fiber sequence (\cite[Theorem 4.12]{Bac22})
\begin{align}
\label{eq:Bac412}
\mathbf{1}[\eta^{-1}]_{(2)}
\to
\mathbf{kw}_{(2)}
\to
\mathbf{kw}_{(2)}[4]
\end{align}
we deduce that $\mathbf{1}[\eta^{-1}]_{(2)}$ is $i_D$-pure.
\endproof



\end{document}